\documentclass{article}
\usepackage{amsfonts,latexsym}
\usepackage{amsmath, times}
\usepackage{amssymb}
\usepackage[colorlinks=true, citecolor=blue, linkcolor=blue, urlcolor=black]{hyperref}
\newcommand{\R}{\mathbb R}
\newcommand{\CC}{\mathbb C}
\newcommand{\N}{\mathbb N}

\newtheorem{theorem}{Theorem}[section]
\newtheorem{lemma}{Lemma}[section]

\newtheorem{definition}{Definition}[section]

\newlength{\defbaselineskip}
\newcommand{\setlinespacing}[2]%
          {\setlength{\baselineskip}{#1 \defbaselineskip}}

\makeatletter

\@addtoreset{equation}{section} \makeatother 
\begin{document}
\begin{center}
 {\Large   {Nodal solutions for  Logarithmic weighted $N$-Laplacian problem with  exponential
nonlinearities}}
\end{center}
\vspace{0.2cm}
\begin{center}
   Brahim Dridi $^{a,b}$ and  Rached Jaidane $^{c}$

 \
\noindent\footnotesize   $^{a}$
Umm Al-Qura University, College of first common year, Department of mathematics, P.O. Box $14035$,
Holly Makkah $21955$, Saudi Arabia\\
 \noindent\footnotesize$^{b}$  University of Tunis El Manar, El Manar preparatory institute for engineering studies, Department of mathematics,  Tunisia.\\
Address e-mail: dridibr@gmail.com\\
\noindent\footnotesize $^c$ Department of Mathematics, Faculty of Science of Tunis, University of Tunis El Manar, Tunisia.\\
 Address e-mail: rachedjaidane@gmail.com\\
\end{center}

\vspace{0.5cm}
\noindent {\bf Abstract.}
In this article, we study the  following problem
$$-\textmd{div} (\omega(x)|\nabla u|^{N-2}  \nabla u) = \lambda\ f(x,u) \quad\mbox{ in }\quad B, \quad u=0 \quad\mbox{ on } \quad\partial B,$$
 where $B$ is the unit ball of $\mathbb{R^{N}}$, $N\geq2$ and  $ w(x)$
a singular weight of logarithm type. The reaction source
$f(x,u)$ is a radial function with respect to $x$ and is subcritical or critical with respect to a maximal growth of  exponential type. By using the constrained minimization in Nehari set coupled with the quantitative deformation lemma and degree theory, we prove the existence of nodal solutions.\\

\noindent {\footnotesize\emph{Keywords:} Weighted Sobolev space, $N$-Laplacian operator, Critical exponential growth, Nodal solutions.\\
\noindent {\bf $2010$ Mathematics Subject classification}: $35$J$20$, $49$J$45$, $35$K$57$, $35$J$60$.}

\section{Introduction and Main results}
In  this paper, we consider  the  following  elliptic problem involving
 logarithmic weighted $N$-Laplacian:
\begin{equation}\nonumber
  \displaystyle (P_{\lambda})~~~~\left\{
      \begin{array}{rclll}
-\textmd{div} (\omega(x)|\nabla u|^{N-2}  \nabla u) &=& \lambda\ f(x,u)& \mbox{in} & B \\
        u&=&0 &\mbox{on }&  \partial B,
      \end{array}
    \right.
\end{equation}

where $B=B(0,1)$ is the unit open ball in $\R^{N}$, $N > 2$, $\lambda$ is a positive parameter, the weight function $w(x)$ is given by
\begin{equation}\label{eq:1.2}
\omega(x)=\Big(\log \frac{1}{|x|}\Big)^{\beta(N-1)}~~\mbox{or}~~~\omega(x)=\Big(\log \frac{e}{|x|}\Big)^{\beta(N-1)}\beta\in[0,1).
\end{equation}
We assume that the nonlinearity $f(x,t):\overline{B}\times\R\rightarrow \R$ is a radial in $x$, continuous function and behaves like $\exp\{\alpha t^{\frac{N}{(N-1)(1-\beta)}}\} ~~\mbox{as}~~t\rightarrow+ \infty ~~$, for some $\alpha >0$ and $\beta\in[0,1)$.\\
Such an equation may arise  in many fields of physics, such
as in non-Newtonian fluids, reaction diffusion problem, turbulent
flows in porous media and image treatment \cite{A,AM,PM,VR}. Here we just give some references which are close to the problem we consider in this note.\\
Without the weight $w(x)$, the problem $(P_{\lambda})$  has been widely studied by several
authors  with different nonlinearities.
G. M. Figueredo and F.B. M. Nunes  in \cite{FN}, consider the following equation
\begin{equation}\label{p}
-\textmd{div}(a(|\nabla u|^{p})|\nabla u|^{p-2}\nabla u)=f(u)~~~ \mbox{ in }~~~ \Omega,~~ u=0 ~~~\mbox{ on } ~~~\partial \Omega,
\end{equation}
where $\Omega\subset\R^{N}$ is bounded, $1<p<N$, the nonlinearity $f:\R\rightarrow \R$ is a superlinear  continuous function with exponential subcritical or exponential critical growth and the function $a$ is $C^{1}$.  By using the  minimization argument and deformation lemma, the authors  proved the existence of a least energy nodal solutions for the equation (\ref{p}) with two nodal domains.\\ Recently, when $a(|\nabla u|^{p})=1$, X. sun and Y. song, see \cite{S}, studied the problem (\ref{p}) in an open smooth bounded domain in the Heisenberg group $\mathbb{H}^{n}=\CC^{n}\times \R$. We also mention the  work of Y. Zhang, Y Yang and  S. Liang, see \cite{ZY}, where they established the existence of changing-sign solutions to the problem (\ref{p}) under logarithmic and exponential nonlinearities.\\
  Weighted $N$-Laplacian elliptic problems of the following type
\begin{equation}\nonumber
\displaystyle \left\{
\begin{array}{rclll}
-\textmd{div} (\omega(x)|\nabla u|^{N-2}  \nabla u) &=&  \displaystyle f(x,u)& \mbox{in} & B \\
u &>&0 &\mbox{in }& B\\
u&=&0 &\mbox{on }&  \partial B,
\end{array}
\right.
\end{equation}
where the weight function $w(x)$ is given in (\ref{eq:1.2}) and  the nonlinearity $f(x,t):\overline{B}\times\R\rightarrow \R$ is positive,  have been investigated in literatures
(see \cite{CRS,DHT,BJ,Z}  and the references therein). We notice that,
the influence of weights on limiting inequalities of Trudinger-Moser type has been studied with some detail in  \cite{C,CR1,CR2,CR3}
 and as consequence the weights have an important impact in the Sobolev norm.\\
Let $\Omega \subset \R^{N}$ be a bounded domain and $w\in L^{1}(\Omega)$ be a nonnegative function, the weighted sobolev space is defined as
$$ W_{0}^{1,N}(\Omega,w)=cl\{u\in
C_{0}^{\infty}(\Omega)/\int_{B}|\nabla u|^{^{N}}w(x)dx <\infty\}$$
For different proprieties and embedding results for the weighted Sobolev spaces, we can refer to \cite{Kuf}.\\
Since the logarithmic weights have a particular significance and considered as the limiting situations of the embedding of the spaces
$W_{0}^{1,N}(\Omega,w)$, we will consider the weight defined by (\ref{eq:1.2}).
We then, restrict our attention to radial functions and consider the subspace
 \begin{equation*}
\mathbf{E}=W_{0,rad}^{1,N}(B,w)=cl\{u \in
C_{0,rad}^{\infty}(B)~~|~~\int_{B}|\nabla u|^{^{N}}\omega(x)dx <\infty\}
 \end{equation*}
 endowed with the norm $$\|u\|=\displaystyle\big(\int_{B}|\nabla u|^{N}\omega(x)dx\big)^{\frac{1}{N}}.$$
 The choice of the weight induced in (\ref{eq:1.2}) and the space $\mathbf{E}$ are also motivated by the following exponential inequalities.
\begin{theorem}\cite{CR2} \label{th1.1}
 Let $\beta\in[0,1)$ and let $w$ be given by (\ref{eq:1.2}), then
  \begin{equation}\label{eq:1.3}
 \int_{B} e^{|u|^{\gamma}} dx <+\infty, ~~\forall~~u\in \mathbf{E},~~
  \mbox{if and only if}~~\gamma\leq \gamma_{N,\beta}=\frac{N}{(N-1)(1-\beta)}=\frac{N'}{1-\beta}
 \end{equation}
and
 \begin{equation}\label{eq:1.4}
 \sup_{\substack{u\in \mathbf{E} \\  \|u\|\leq 1}}
 \int_{B}~e^{\alpha|u|^{\gamma_{N,\beta} }}dx < +\infty~~~~\Leftrightarrow~~~~ \alpha\leq \alpha_{N,\beta}=N[\omega^{\frac{1}{N-1}}_{N-1}(1-\beta)]^{\frac{1}{1-\beta}},
 \end{equation}
where $\omega_{N-1}$ is the area of the unit sphere $S^{N-1}$ in $\R^{N}$ and $N'$ is the H$\ddot{o}$lder conjugate of $N$.
 \end{theorem}

Let $\gamma:= \gamma_{N,\beta}=\displaystyle\frac{N'}{1-\beta}$,  in view  of inequalities (\ref{eq:1.5}) and (\ref{eq:1.6}), we say that $f$ has subcritical growth at $+\infty$ if
\begin{equation} \label{eq:1.5}
\lim_{s\rightarrow +\infty}\frac{|f(x,s)|}{e^{\alpha s^{\gamma}}}=0,~~~~\mbox{for all}~~ \alpha>0
\end{equation}
and $f$ has critical growth at $+\infty$ if there exists some $\alpha_{0}>0$,
\begin{equation}\label{eq:1.6}
\lim_{s\rightarrow +\infty}\frac{|f(x,s)|}{e^{\alpha s^{\gamma}}}=0,~~~\forall~\alpha>\alpha_{0}
~~\mbox{and}~~~~\lim_{s\rightarrow +\infty}\frac{|f(x,s)|}{e^{\alpha s^{\gamma}}}=+\infty,~~\forall~\alpha>\alpha_{0}.\\
\end{equation}

In this paper, we deal  with problem $(P_{\lambda})$ under subcritical and critical growth nonlinearities.
Furthermore, we suppose that $f(x,t)$ satisfies the following hypothesis:
\begin{enumerate}
\item[$(V_{1})$] $f: B \times \mathbb{R}\rightarrow\mathbb{R}$ is $C^{1}$ and radial in $x$.
\item[($V_{2})$] There exist $\theta > N$  such that  we have
 $$0 < \theta F(x,t)\leq tf(x,t), \forall (x, t)\in~~B\times\mathbb{R} \setminus\{0\} $$
  where
$$F(x,t)=\displaystyle\int_{0}^{t}f(x,s)ds.$$
\item [$(V_{3})$]  For each $x\in B$,~$\displaystyle t \mapsto \frac{f(x,t)}{|t|^{N-1}}~~\mbox{is increasing for}~~ t\in~~\mathbb{R} \setminus\{0\}$.
 \item [$(V_{4})$]   $\displaystyle\lim_{t\rightarrow 0}\frac{|f(x,t)|}{|t|^{N-1}}=0.$
\end{enumerate}
We give an example of such nonlinearity. The nonlinearity $f(x,t)=|t|^{N-1}t+|t|^{N}t\exp(\alpha |t|^{\gamma})$ satisfies the assumptions $(V_{1})$, $(V_{2})$, $(V_{3})$ and $(V_{4})$ .\\


We will consider the following definition of  solutions.\\
\begin{definition}\label {def1.1}
We say that a function $u\in \mathbf{E}$ is a weak solution of the problem $(P_{\lambda})$ if
\begin{equation*}\label {eq:1.9}
\int_{B}|\nabla u|^{N-2}\nabla u .\nabla \varphi~ w(x) dx =\lambda
\int_{B}f(x,u) \varphi dx,~~\forall~\varphi \in \mathbf{E}.
\end{equation*}
\end{definition}
Let $\mathcal{J_{\lambda}} :\mathbf{E} \rightarrow \R$ be the functional given by
 \begin{equation}\label{eq:1.7}
\mathcal{J_{\lambda}}(u)=\frac{1}{N}\int_{B}|\nabla u|^{N}w(x)dx-\lambda\int_{B}F(x,u)dx,
\end{equation}
where
$$ F(x,t)=\displaystyle\int_{0}^{t}f(x,s)ds.$$
The energy functional $\mathcal{J_{\lambda}}$ is well defined and of class $C^{1}$ since there exist $ a,~ C>0$ positive constants and there exists $t_{1} >1$ such for that
\begin{equation*}\label {eq:1.11}
|f(x,t)|\leq C e^{a ~t^{\gamma}}, ~~~~~~\forall |t| >t_{1},
\end{equation*}
whenever the nonlinearity $f(x,t)$ is critical or subcritical at $+\infty$.\\It is quite clear that finding non trivial weak solutions to the problem $(P_{\lambda})$ is equivalent to finding non-zero critical points of the functional $\mathcal{J_{\lambda}}$. Moreover, we have $$\langle\mathcal{J}'_{\lambda}(u),\varphi\rangle=\mathcal{J_{\lambda}'}(u)\varphi=\int_{B}\big(\omega(x)~|\nabla u|^{N-2}\nabla u~\nabla \varphi\big)~dx-\lambda\int_{B}f(x,u)~ \varphi~dx~,~~\varphi \in\mathbf{E}\cdot$$ \\
We define the Nehari set as
$$\displaystyle\ \mathcal{N}_{\lambda}:=\{u\in \mathbf{E}:\langle \mathcal{J}'_{\lambda}(u),u^{+}\rangle=\langle\mathcal{J}'_{\lambda}(u),u^{-}\rangle=0, u^{+}\neq 0, u^{-}\neq 0 \},$$
where $u^{+}=\max \{u(x), 0\}$,  $u^{-}=\min \{u(x), 0\}$.\\It's easy to verify the following decomposition $$\mathcal{J}_{\lambda}(u)=\mathcal{J}_{\lambda}(u^{+})+\mathcal{J}_{\lambda}(u^{-}),  $$ and
$$\langle\mathcal{J}'_{\lambda}(u),u^{+}\rangle=\langle\mathcal{J}'_{\lambda}(u^{+}),u^{+}\rangle~~ \mbox{and}~~\langle\mathcal{J}'_{\lambda}(u),u^{-}\rangle=\langle\mathcal{J}'_{\lambda}(u^{-}),u^{-}\rangle$$
We also give the following definitions of the so called nodal solutions and least energy sign-changing solution of problem $(P_{\lambda})$.
\begin{definition}\begin{itemize}
\item $\upsilon\in \mathbf{E}$ is called nodal or sign-changing solution of problem $(P_{\lambda})$ if $\upsilon$ is a solution of problem $(P_{\lambda})$ and $\upsilon^{\pm}\neq  0~~\mbox{a.e in}~~ B$.\item $\upsilon\in \mathbf{E}$
is called least energy sign-changing solution of problem $(P_{\lambda})$ if $\upsilon$  is a sign-changing solution of $(P_{\lambda})$ and
\begin{equation*}
\mathcal{J_{\lambda}}(\upsilon)=\displaystyle\inf\{\mathcal{J_{\lambda}}(u):\mathcal{J'_{\lambda}}(u)=0, u^{\pm}\neq0~~\mbox{a.e in}~~B \}
\end{equation*}
\end{itemize}
\end{definition}
Influenced by the works cited above, we try to find a minimize of the energy functional $\mathcal{J_{\lambda}}$ over the following minimization problem,
\begin{equation*}\label {eq:1.12}
\displaystyle c_{\lambda}=\inf_{u\in\mathcal{N}_{\lambda}}\mathcal{J_{\lambda}}(u)
\end{equation*}
Our approach is based on the
Nehari manifold method, which was introduced in \cite{Kav} and is by now a well-
established and useful tool in finding solutions of problems with a variational
structure, see \cite{FMR}.\\
To our best knowledge, there are few results for the nodal solutions to the N-weighted Lapalace equation with critical exponential nonlinearity on the weighted Sobolev space $\mathbf{E}$.\\
Now, we give our main results as follows:
\begin{theorem}\label{th1.2}~~ Let $f(x,t)$ be a function that has a subcritical growth at
$+\infty$   $(V_{1})$, $(V_{2})$, $(V_{3})$, and $(V_{4})$ are satisfied.
For $\lambda>0$,  the  problem $(P_{\lambda})$ has a least energy nodal (sign-changing) radial solution $\upsilon\in\mathcal{N}_{\lambda}$ .
 \end{theorem}
For a critical growth nonlinearity, the following result holds.
\begin{theorem}\label{th1.3}~~
Assume that $f(x,t)$ has a  critical growth at $+\infty$ for some $\alpha_{0}$ and  $(V_{1})$, $(V_{2})$, $(V_{3})$ and $(V_{4})$ are satisfied. Then, there exists $\lambda^{*}>0$ such that for $\lambda>\lambda^{*}$,
problem $(P_{\lambda})$ has a least energy nodal (sign-changing) radial solution $\upsilon\in\mathcal{N}_{\lambda}$.
\end{theorem}

This present work is organized as follows: in section $2$, some preliminaries for the compactness analysis
 are presented. In section $3$, we give  some technical key lemmas . In section $4$ we prove our result in the subcritical case. Section $5$ is devoted for the critical case which is more difficult. We use a concentration compactness result of Lions type to prove Theorem \ref{th1.3}.\\
Finally,  we note that a constant $C$ may change from line to another and sometimes we index the constants in order to show how they change.


\section{  Preliminaries for the compactness analysis}
In this section, we will present a number of technical Lemmas for our future use. We begin by the radial
Lemma.

\begin{lemma}\cite{CR2}\label{lemr}Let $u$ be a radially symmetric
 function in $C_{0}^{1}(B)$. Then, we have
 $$|u(x)|\leq \displaystyle\frac{|\log(|x|)|^{\frac{{1-\beta}}{N'}}}{\omega_{N-1}^{\frac{1}{N}}(1-\beta)^{\frac{1}{N'}}}\|u\|,$$
where $\omega_{N-1}$ is the area of the unit sphere $S^{N-1}$ in $\mathbb{R}^{N}$.
\end{lemma}
 It follows that the embedding  $\mathbf{E}\hookrightarrow L^{q}(B)$ is continuous  for all $q \geq 1$, and that there exists a constant $C>0$ such that $\|u\|_{N'q}\leq C \|u\|$, for all $u\in \mathbf{E}$. Moreover, the embedding  $\mathbf{E}\hookrightarrow L^{q}(B)$ is compact for all  $q \geq N$.
\begin{lemma}\label{Lionstype}  Let
$(u_{k})_{k}$  be a sequence in $\mathbf{E}$. Suppose that,  $\|u_{k}\|=1$, $u_{k}\rightharpoonup u$ weakly in $\mathbf{E}$, $u_{k}(x)\rightarrow u(x) ~~a.e.$ in $ B$, $\nabla u_{k}(x)\rightarrow \nabla u(x) ~~a.e.$ in $B$ and $u\not\equiv 0$. Then
$$\displaystyle\sup_{k}\int_{B}e^{p~\alpha_{N,\beta}
|u_{k}|^{\gamma}}dx< +\infty,~~\mbox{where}~~ \alpha_{N,\beta}=N[w^{\frac{1}{N-1}}_{N-1}(1-\beta)]^{\frac{1}{1-\beta}},$$
for all $1<p< \mathcal{U}(u)$,  where $\mathcal{U}(u)$ is given by:
 $$\mathcal{U}(u):=\displaystyle \left\{
      \begin{array}{rcll}
&\displaystyle\frac{1}{(1-\|u\|^{N})^{\frac{\gamma}{N}}}& \mbox{ if }\|u\| <1\\
       &+\infty& \mbox{ if } \|u\|=1.\\
 \end{array}
    \right.$$
\end{lemma}
Proof: For $a, b \in \R, q>1$. If $q'$ its conjugate i.e. $\frac{1}{q}+\frac{1}{q'}=1$
we have, by young inequality, that
$$e^{a+b}\leq \frac{1}{q}e^{qa}+ \frac{1}{q'}e^{q'b}.$$

Also, we have
\begin{equation*}\label {eq:3.1}
(1+a)^{q}\leq (1+\varepsilon) a^{q}+(1-\frac{1}{(1+\varepsilon)^{\frac{1}{q-1}}})^{1-q},~~\forall a\geq0,~~\forall\varepsilon>0~~\mbox{and}~~\forall q>1.
\end{equation*}
So, we get
$$
      \begin{array}{rcll}
|u_{k}|^{\gamma}&=& |u_{k}-u+u|^{\gamma}\\
&\leq& (|u_{k}-u|+|u|)^{\gamma}\\
&\leq& (1+\varepsilon)|u_{k}-u|^{\gamma}+\big(1-\frac{1}{(1+\varepsilon)^{\frac{1}{\gamma-1}}}\big)^{1-\gamma}|u|^{\gamma}\\
 \end{array}
   $$
which implies that
 \begin{align*}
\int_{B}e^{p~\alpha_{N,\beta}
|u_{k}|^{\gamma}}dx &\leq
\frac{1}{q}\int_{B} e^{pq~\alpha_{N,\beta}
(1+\varepsilon)|u_{k}-u|^{\gamma}}dx\\
&+\displaystyle\frac{1}{q'}\int_{B} e^{pq'~\alpha_{N,\beta}
(1-\frac{1}{(1+\varepsilon)^{\frac{1}{\gamma-1}}})^{1-\gamma}|u|^{\gamma}}dx,
\end{align*}
for any $p>1$. Since $(1-\frac{1}{(1+\varepsilon)^{\frac{1}{\gamma-1}}})^{1-\gamma}\leq1$,
then$$
\displaystyle\frac{1}{q'}\int_{B} e^{pq'~\alpha_{N,\beta}
(1-\frac{1}{(1+\varepsilon)^{\frac{1}{\gamma-1}}})^{1-\gamma}|u|^{\gamma}}dx\leq
\displaystyle\frac{1}{q'}\int_{B} e^{pq'~\alpha_{N,\beta}
|u|^{\gamma}}dx=
\displaystyle\frac{1}{q'}\int_{B} e^{\big((pq'~\alpha_{N,\beta})^{\frac{1}{\gamma}}
|u|\big)^{\gamma}}dx.$$\\
From (\ref{eq:1.5}), the last integral is finite.\\ To complete the
evidence, we have to prove that for every $p$ such that $1<p<\mathcal{U}(u)$,
\begin{equation}\label {eq:2.1}
\sup_{k}\int_{B} e^{pq~\alpha_{N,\beta}
(1+\varepsilon)|u_{k}-u|^{\gamma}}dx<+\infty,
\end{equation}
 for some $\varepsilon>0$ and $q>1$.\\
In the following, we suppose that
$\|u\|<1$ and in the case of $\|u\|=1$, the proof  is similar.\\
When $$\|u\|<1$$
and
$$p<\displaystyle\frac{1}{(1-\|u\|^{N})^{\frac{\gamma}{N}}},$$
there exists $\nu>0$ such that
$$p(1-\|u\|^{N})^{\frac{\gamma}{N}}(1+\nu)<1.$$
On the other hand,
 From the Brezis-Lieb's lemma \cite{Br} it holds that
 \begin{equation*}\label {eq:2.3}
 \|u_{k}-u\|^{N}=\|u_{k}\|^{N}-\|u\|^{N}+o(1)~~\mbox{where}~~
  o(1)\rightarrow 0~~ \mbox{as}~~k\rightarrow +\infty.
\end{equation*}
   Then,
   $$\|u_{k}-u\|^{N}=1-\|u\|^{N}+o(1),$$
so,
$$\displaystyle\lim_{k\rightarrow+\infty}\|u_{k}-u\|^{\gamma}=(1-\|u\|^{N})^{\frac{\gamma}{N}}.$$

 Therefore, for every
$\varepsilon>0$, there exists $k_{\varepsilon}\geq 1$ such that
$$\|u_{k}-u\| ^{\gamma}\leq (1+\varepsilon)(1-\|u\|^{N})^{\frac{\gamma}{N}},~~\forall ~~k\geq k_{\varepsilon}.$$
Then, for $q=1+\varepsilon$ with
$\varepsilon=\sqrt[3]{1+\nu}-1$ and for any  $ k\geq k_{\varepsilon}$, we have
 $$pq(1+\varepsilon)\|u_{k}-u\|^{\gamma}\leq 1.$$
Consequently,
$$\begin{array}{rlll}
\displaystyle\int_{B} e^{pq~\alpha_{N,\beta}
(1+\varepsilon)|u_{k}-u|^{\gamma}}dx&\leq&
\displaystyle\int_{B} e^{
(1+\varepsilon)pq~\alpha_{N,\beta}(\frac{|u_{k}-u|}{\|u_{k}-u\|})^{\gamma}\|u_{k}-u\|^{\gamma}}dx\\
&\leq&\displaystyle\int_{B} e^{~\alpha_{N,\beta}(\frac{|u_{k}-u|}{\|u_{k}-u\|})^{\gamma}}dx\\
&\leq &\displaystyle\sup _{\|u\|\leq 1}\displaystyle\int_{B}
e^{~\alpha_{N,\beta}|u|^{\gamma}}dx <+\infty.
\end{array}
$$
Now, (\ref{eq:2.1})  follows from  (\ref{eq:1.4}). This complete the proof.
A second important Lemma.
 \begin{lemma}\cite{FMR}\label{lem2}
 Let $\Omega\subset \mathbb{R^{N}}$ be a bounded domain and $f:\overline{\Omega}\times\mathbb{R}$
  be a continuous function. Let $\{u_{n}\}_{n}$ be a sequence in $L^{1}(\Omega)$
converging to $u$ in $L^{1}(\Omega)$. Assume that $\displaystyle f(x,u_{n})$ and
$\displaystyle f(x,u)$ are also in $ L^{1}(\Omega)$. If
$$\displaystyle\int_{\Omega}|f(x,u_{n})u_{n}|dx \leq C,$$ where $C$
is a positive constant, then $$f(x,u_{n})\rightarrow
f(x,u)~~in~~L^{1}(\Omega).$$
\end{lemma}

\section{ Some technical lemmas }
In the following we assume, unless otherwise stated, that the function $f$ satisfies the conditions $(V_{1})$ to $(V_{4})$. Let $u$ $\mathbf{E}$ with $u^{\pm}\not\equiv 0 $ a.e. in the ball $ B$, and we deﬁne the function
$\displaystyle\Upsilon_{u} : [0, \infty) \times [0, \infty) \rightarrow\R$
 and mapping
$L_{u} : [0, \infty) \times [0, \infty) \rightarrow\R^{2} $ as
\begin{equation}\label{eq:3.1}
 \displaystyle\Upsilon_{u}(p, q) = \mathcal{J}_{\lambda}(pu^{+} + qu^{-}),
\end{equation}
and
\begin{equation}\label{eq:3.2}
 \displaystyle L_{u}(p, q) = \left (\langle \mathcal{J'}_{\lambda}(pu^{+} + qu^{-}),  pu^{+} \rangle ,  \langle \mathcal{J'}_{\lambda} (pu^{+} + qu^{-}, qu^{-}\rangle)  \right.	
\end{equation}
\begin{lemma}\label{lem1}\begin{itemize}
\item [(i)]For  each  $u$ $\mathbf{E}$ with  $u^{+}\neq0$  and  $u^{-}\neq 0$,  there  exists  an  unique couple
 $\displaystyle(p_{u}, q_{u})\in (0,\infty) \times(0,\infty)$ such that
$\displaystyle p_{u}u^{+} + q_{u}u^{-} \in\mathcal{N}_{\lambda}.$ In particular,  the set $\mathcal{N}_{\lambda}$ is nonempty.
\item[(ii)]  For all $p, q \geq 0 $ with $(p, q) \neq (p_{u}, q_{u}),$ we have
$$\displaystyle \mathcal{J}_{\lambda}(pu^{+} + qu^{-}) <  \mathcal{J}_{\lambda}(p_{u}u^{+} + q_{u}u^{-})\cdot$$\end{itemize}

\end{lemma}
Proof.$(i)$ \\
Since $f$ is subcritial or critical, and 	From $(V_{1})$ and $(V_{4})$,  for all $\varepsilon> 0$, there exists a positive constant $C_{1} = C_{1}(\epsilon )$ such that
\begin{equation}\label{eq:3.3}
f(x,t)t\leq \varepsilon |t|^{N} +C_{1}|t|^{s }\exp(\alpha|t|^{\gamma}  )\mbox{ for all }	\alpha > \alpha_{0}, s > N.
\end{equation}
 Now, given $u\in\mathbf{E}$ fixed with  $u^{+}\neq0$  and  $u^{-}\neq 0$. From (\ref{eq:3.3}), for all $\varepsilon>0$, we have
\begin{align}\label{eq:3.4}
 \displaystyle\langle \mathcal{J}'_{\lambda}(pu^{+} + qu^{-}),  pu^{+} \rangle &=\langle \mathcal{J}_{\lambda}(pu^{+} ),  pu^{+} \rangle \nonumber  \\
 & =\|pu^{+}\|^{N}-\lambda\int_{B}f(x,pu^{+})pu^{+}dx\nonumber   \\
 & \geq \|pu^{+}\|^{N} \nonumber
 -\lambda\epsilon \int_{B}|pu^{+}|^{N}dx\nonumber -\lambda C_{1}\int_{B}|pu^{+}|^{s}\exp (\alpha p|u^{+}|^{\gamma})dx
 \end{align}
   Using the
 H\"{o}lder inequality, with $a, a' > 1$ such that $\displaystyle\frac{1}{a} + \frac{1 }{a'}  = 1$, and Lemma \ref{lemr}, we get
   \begin{align}
  \displaystyle\langle \mathcal{J}'_{\lambda}(pu^{+} + qu^{-}),  pu^{+} \rangle  & \geq \|pu^{+}\|^{N}-\lambda\epsilon C_{2}\|pu^{+}\|^{N}\nonumber
   -\lambda C_{1}\left( \int_{B}|pu^{+}|^{a's}dx\right)^{ \frac{1 }{a'} }\left(\int_{B}\exp (\alpha p a|u^{+}|^{\gamma})dx\right)^{\frac{1}{a}}\nonumber \\ &\geq \left(1-\epsilon C_{2}-\lambda \epsilon C_{1} \right)\|pu^{+}\|^{N}-\lambda C_{1}\left(\int_{B}\exp \big(\alpha  a\|pu^{+}\|^{\gamma}\big(\frac{|u^{+}|}{\|u^{+}\|}\big)^{\gamma}\big)dx\right)^{\frac{1}{a}} C_{3}\|pu^{+}\|^{s}\nonumber
   \end{align}
  By (\ref{eq:1.4}), the last integral is finte provided $p>0$ is chosen small enough such that  $\displaystyle\alpha a\|pu^{+}\|^{\gamma}\leq \alpha_{N,\beta}$. Then,
\begin{align}
  \displaystyle\langle \mathcal{J}'_{\lambda}(pu^{+} + qu^{-}),  pu^{+} \rangle  & \geq \left(1-\epsilon C_{2}-\lambda \epsilon C_{1} \right)\|pu^{+}\|^{N}-\lambda C_{4}\|pu^{+}\|^{s}\end{align}
holds. Choosing  $\epsilon > 0$ such that $1-\epsilon C_{2}-\lambda \epsilon C_{1} > 0$ anf for small  $p > 0 $ and for all  $q > 0$ and  $ s > N$,  we get
$\displaystyle\langle \mathcal{J}'_{\lambda}(pu^{+} + qu^{-}),  pu^{+} \rangle > 0$ . In the similar way, it can be proved that
 $\displaystyle\langle \mathcal{J}'_{\lambda}(pu^{+} + qu^{-}),  pu^{-} \rangle > 0$ for $q > 0$ small enough and all $p > 0$.
Therefore, it is quite easy to state that there exists $t_{1} > 0$ such that
\begin{equation}\label{eq:3.5}
\displaystyle\langle \mathcal{J}'_{\lambda}(t_{1}u^{+} + qu^{-}),  t_{1}u^{+} \rangle >0,
\quad \langle \mathcal{J}'_{\lambda}(pu^{+} + t_{1}u^{-}),  t_{1}u^{-} \rangle >0\\  \mbox{  for all  } p, q > 0.
\end{equation}
 From $(V_{3})$, we can derive that there exists $C_{5}, C_{6} > 0 $ such that
\begin{equation}\label{eq:3.6}
F (x,t)\geq  C_{5}|t|^{\theta}-C_{6}.
\end{equation}
Now, choose $ p = t_{2}^{\ast}  > t_{1}$  with $t_{2}^{\ast}$  large enough. Then, by using $(\ref{eq:3.3})$, $(\ref{eq:3.6})$, we get
\begin{align*}
 \displaystyle\langle \mathcal{J}'_{\lambda}(t_{2}^{\ast}u^{+} + qu^{-}),  t_{2}^{\ast}u^{+} \rangle
  &=\langle \mathcal{J'}_{\lambda}(t_{2}^{\ast}u^{+} ),  t_{2}^{\ast}u^{+} \rangle \nonumber  \\
  & \leq \|t_{2}^{\ast}u^{+}\|^{N}
  -\lambda\int_{B}C_{5}|t_{2}^{\ast}u^{+}|^{\theta}dx+\lambda C_{6}|B|\nonumber  \\
  &\leq 0,
\end{align*}
for $q \in [t_{1}, t_{2}^{\ast}]$. Also, we can choose $q = t_{2}^{\ast} >t_{1}$  with $t_{2}^{\ast}$  large enough and then
$$ \displaystyle \langle\mathcal{J}'_{\lambda}(t_{2}^{\ast}u^{+} +t_{2}^{\ast}u^{-}),  t_{2}^{\ast}u^{+} \rangle <0 \mbox{   holds for }   p\in [t_{1}, t_{2}^{\ast}].$$

Therefore, if $t_{2} > t_{2}^{\ast}$ is large enough, then we obtain that
\begin{equation}\label{eq:3.7}
 \displaystyle\mathcal{J}'_{\lambda}(t_{2}u^{+} +qu^{-}),  t_{2}u^{+} \rangle <0 \quad\mbox{ and }\quad \langle \mathcal{J}'_{\lambda}(pu^{+} +t_{2}u^{-}),  t_{2}u^{-} \rangle <0
\mbox{ for all }   p, q\in [t_{1}, t_{2}].
\end{equation}
Joining $(\ref{eq:3.5} )$ and $(\ref{eq:3.7})$ with Miranda's  Theorem \cite{Av}, there exists at least a couple of points $(p_{u}, q_{u})\in (0, \infty)\times (0, \infty)$
 such that $L_{u}(p_{u}, q_{u})=(0, 0)$ , i.e, $p_{u}u^{+} +q_{u}u^{-} \in\mathcal{ N}_{\lambda}.$\\
Now we will show the uniqueness of the couple $(p_{u}, q_{u})$. Indeed, it is suﬃcient to show that if
 $u \in\mathcal{ N}_{\lambda}$ and  $p_{0}u^{+} +q_{0}u^{-} \in\mathcal{ N}_{\lambda}$ with $ p_{0 } >  0$  and  $q_{0 } >  0$,  then  $(p_{0}, q_{0})  =  (1, 1)$.
 Let us assume that $u \in\mathcal{ N}_{\lambda}$ and  $p_{0}u^{+} +q_{0}u^{-} \in\mathcal{ N}_{\lambda}.$We will then get
 $\displaystyle\langle\mathcal{J}'_{\lambda}(p_{0}u^{+} +q_{0}u^{-}),  p_{0}u^{+} \rangle= 0,$
 $\displaystyle\langle\mathcal{J}'_{\lambda}(p_{0}u^{+} +q_{0}u^{-}),  p_{0}u^{-} \rangle= 0,$ and
$\displaystyle\langle\mathcal{J}'_{\lambda}(u), u^{\pm} \rangle= 0,$ that is,

\begin{equation}\label{eq:3.8}
 \displaystyle \|p_{0}u^{+}\|^{N} =\lambda\int_{B}f(x,p_{0}u^{+})p_{0}u^{+}dx
\end{equation}
\begin{equation} \label{eq:3.9}
  \|b_{0}u^{-}\|^{N}=\lambda\int_{B}f(x,q_{0}u^{-})q_{0}u^{-}dx
\end{equation}
\begin{equation} \label{eq:3.10}
  \|u^{+}\|^{N}=\lambda\int_{B}f(x,u^{+})u^{+}dx
\end{equation}
\begin{equation}\label{eq:3.11}
   \|u^{-}\|^{N}=\lambda\int_{B}f(x,u^{-})u^{-}dx
\end{equation}
Combining (\ref{eq:3.8}) and (\ref{eq:3.10}), we deduce that
 $$ 0=\lambda\int_{B}\frac{f(x,p_{0}u^{+})p_{0}u^{+}}{p_{0}^{N}}dx-\lambda\int_{B}f(x,u^{+})u^{+}dx.$$
 It follows from $(V_{4})$  that $t\mapsto \frac{f(x,t)}{t^{N-1}}$ is increasing for $t > 0$, which implies
that $p_{0} = 1$. We can also show, using $(V_{4})$, (\ref{eq:3.9}) and (\ref{eq:3.11}), that $q_{0} = 1$. This completes the proof of $(i)$.\\
$(ii)$ To prove $(ii)$, it is sufficient to show that $(p_{u}, q_{u})$ is the unique maximum point of $\displaystyle\Upsilon_{u} \in [0, \infty) \times [0, \infty) $. From (\ref{eq:3.7}), (\ref{eq:3.8}) and $\theta > N $, we have
\begin{align*}
  \displaystyle\Upsilon_{u}(p,q) &= \mathcal{J}_{\lambda}(pu^{+} +qu^{-}) \\
  &= \frac{1}{N} \|pu^{+} +qu^{-}\|^{N}-\lambda\int_{B}F(x,pu^{+} +qu^{-})dx \\
  &\leq \frac{p^{N}}{N} \|u^{+}\|^{N}+\frac{q^{N}}{N} \|u^{-}\|^{N}-\lambda C_{5}p^{\theta}\int_{B}|u^{+}|^{\theta}dx-
  \lambda C_{5}q^{\theta}\int_{B}|u^{-}|^{\theta}dx+\lambda C_{6}|B|
\end{align*}
which implies that  $\displaystyle\lim_{|(p,q)|\rightarrow \infty} \Upsilon_{u}(p,q)=-\infty$.
Hence, it suffices to see that the maximum point of $\Upsilon_{u}$ cannot be realized on the boundary of $[0, \infty ) \times[0, \infty).$
We argue by contradiction and assume that  $(0, q)$ with $q\geq 0$ is a maximum point of $\Upsilon_{u}$. Then from (\ref{eq:3.5}), we have
 $$\displaystyle p\frac{d}{dp}[\mathcal{J}_{\lambda}(pu^{+} +qu^{-})]=\langle\mathcal{J}'_{\lambda}(pu^{+}), pu^{+} \rangle >0,$$
 for small $p > 0$, which means that $\Upsilon_{u}$ is increasing with respect to $p$ if $p > 0$ is small enough. This gives a contradiction. We can similarly deduce that $\Upsilon_{u}$ can not realize its global maximum on $(p, 0)$ with $p \geq 0.$
\begin{lemma}\label{lem2}
 For any $u$ $\mathbf{E}$ with  $u^{+}\neq0$  and  $u^{-}\neq 0$, such that
   $\langle\mathcal{J}'_{\lambda}(pu^{+}, pu^{+} \rangle \leq0$, the unique maximum point $(p_{u}, q_{u})$
   of $\Upsilon_{u}$ on $[0, \infty)\times [0, \infty)$ belongs to $(0,1]\times(0,1]$.
\end{lemma}
Proof.
Here we will only prove that $0 < p_{u} \leq 1$. The proof of $0 < q_{u} \leq 1$ is similar. Since $p_{u}u^{+}, q_{u}u^{-}\in\mathcal{N}_{\lambda},$ we have that
\begin{equation}\label{eq:3.12}
  \displaystyle\|p_{u}u^{+}\|^{N}=\lambda\int_{B}f(x,p_{u}u^{+})u^{+}dx
\end{equation}
Moreover, by  $\displaystyle\langle\mathcal{J}'_{\lambda}(pu^{+}, pu^{+} \rangle \leq0$, we have that
\begin{equation}\label{eq:3.13}
    \|u^{+}\|^{N}\leq\lambda\int_{B}f(x,u^{+})u^{+}dx.
\end{equation}
Combining (\ref{eq:3.12}) and (\ref{eq:3.13}), it follows that
\begin{equation}\label{eq:3.14}
\displaystyle
\int_{B}f(x,u^{+})u^{+}dx\geq \int_{B}\frac{f(x,p_{u}u^{+})p_{u}u^{+}}{p_{u}^{N}}dx.
\end{equation}
Now, we suppose, by contradiction, that $p_{u} > 1$. By $(V_{3})$,  $t\mapsto \displaystyle\frac{f(x,t)}{t^{N-1}}$ is
increasing for $t > 0$, which contradicts inequality (\ref{eq:3.14}).Therefore, $0 < p_{u} \leq 1.$	
\begin{lemma}\label{lem3}
 For all $u\in \mathcal{N}_{\lambda}$,
 \begin{itemize}
   \item [$(i)$] there exists $\kappa>0$ such that\\
   $ \|u^{+}\|, \|u^{-}\| \geq \kappa ;$
   \item[$(ii)$] $\mathcal{J}_{\lambda}(u) \geq (\frac{1}{N}-\frac{1}{\theta})\|u\|^{N}$
 \end{itemize}
 \end{lemma}
  Proof. $(i)$
We argue by contradiction.  Suppose that there exists a sequence $\{u_{n}^{+}\} \subset \mathcal{N}_{\lambda} $
such that $u_{n}^{+}\rightarrow 0$ in $ \mathbf{E}.$
Since $\{u_{n}\} \subset \mathcal{N}_{\lambda}$, then  $\displaystyle\langle\mathcal{J}'_{\lambda}(u_{n}) ,u_{n}^{+} \rangle =0$.
Hence, it follows from (\ref{eq:3.3}), (\ref{eq:3.4}) and the radial Lemma \ref{lemr} that

\begin{align}\label{eq:3.15}
\|u_{n}^{+}\|^{N}&=\lambda\int_{B}f(u_{n}^{+})u_{n}^{+}dx \\
 & \leq \epsilon \lambda \int_{B}|u_{n}^{+}|^{N}dx + \lambda C_{1}\int_{B}|u_{n}^{+}|^{s}\exp(\alpha |u_{n}^{+}|^{\gamma} )dx\nonumber\\
 & \leq \epsilon \lambda C_{6} \|u_{n}^{+}\|^{N} + \lambda C_{1}\int_{B}|u_{n}^{+}|^{s}\exp(\alpha |u_{n}^{+}|^{\gamma} )dx\nonumber
\end{align}
Let $a>1$ with $\frac{1}{a}+\frac{1}{a'}=1$. Since $u_{n}^{+}\rightarrow 0\mbox{ in }~~\mathbf{E}$,
for $n$ large enough, we get
$\displaystyle\|u_{n}^{+}\|\leq(\frac{\alpha_{ N,\beta}}{\alpha a})^{\frac{1}{\gamma}}$. From H\"{o}lder inequality, (\ref{eq:1.4}) and again the radial Lemma \ref{lemr},
 we have
\begin{align*}
  \int_{B}|u_{n}^{+}|^{s}\exp(\alpha |u_{n}^{+}|^{\gamma} )dx&\leq
  \left( \int_{B}|u_{n}^{+}|^{sa'}dx\right)^{\frac{1}{a'}}
  \left(\int_{B}\exp \big(\alpha  a\|u^{+}\|^{\gamma}\big(\frac{|u^{+}|}{\|u^{+}\|}\big)^{\gamma}\big)dx\right)^{\frac{1}{a}} \\
  &\leq C_{7} \left( \int_{B}|u_{n}^{+}|^{sa'}dx\right)^{\frac{1}{a'}}\leq C_{8} \|u_{n}^{+}\|^{s}
\end{align*}
Combining (\ref{eq:3.15}) with the last inequality, for $n$ large enough, we obtain

\begin{equation}\label{eq:3.16}
 \|u_{n}^{+}\|^{N}\leq\lambda\epsilon C_{6} \|u_{n}^{+}\|^{N}
 +\lambda C_{8} \|u_{n}^{+}\|^{s}
\end{equation}
Choose suitable $\epsilon > 0 $ such that  $1-\lambda\epsilon C_{6} > 0$.
Since  $N < s$, then (\ref{eq:3.16}) contradicts the fact that $u_{n}^{+}\rightarrow 0\mbox{ in }\mathbf{E}$.\\
$(ii)$ Given $u \in \mathcal{N}_{\lambda}$, by the deﬁnition of $\mathcal{N}_{\lambda}$ and $(V_{3})$ we obtain
\begin{align*}
  \mathcal{J}_{\lambda}(u) &= \mathcal{J}_{\lambda}(u)-\frac{1}{\theta}\langle\mathcal{J}'_{\lambda}(u), u\rangle \\
  &= \frac{1}{N}\|u_{n}\|^{N}
  +\lambda \big(\int_{B}\frac{1}{\theta}f(x,u)u-F(x,u)dx\big)\\
  &\geq (\frac{1}{N}-\frac{1}{\theta})\|u\|^{N}
\end{align*}
Lemma \ref{lem3} implies that $\mathcal{J}_{\lambda}(u) > 0$  for all $u\in \mathcal{N}_{\lambda}$.
As a consequence, $\mathcal{J}_{\lambda}$ is bounded by below in $\mathcal{N}_{\lambda}$, and therefore
$\displaystyle c_{\lambda}:=\inf_{u\in \mathcal{N}_{\lambda}} \mathcal{J}_{\lambda}(u) $
 is well-deﬁned.\\ The following lemma deals with the asymptotic property of $c_{\lambda}$.
\begin{lemma}\label{lem4}
Let $\displaystyle c_{\lambda}=\inf_{u\in \mathcal{N}_{\lambda}} \mathcal{J}_{\lambda}(u)$, then $\displaystyle\lim_{\lambda\rightarrow \infty}c_{\lambda}=0$
\end{lemma}
Proof. Let us Fix $u ~\mathbf{E}$  with $u^{\pm}=0$. Then, by Lemma \ref{lem1}, there exists a point pair $(p_{\lambda}, q_{\lambda})$ such that
$p_{\lambda}u^{+}+q_{\lambda}u_{-} \in \mathcal{N}_{\lambda}$  for each $\lambda > 0$. Let $\mathcal{T}_{u}$ be the set defined by
$$\displaystyle \mathcal{T}_{u} := \{(p_{\lambda}, q_{\lambda})\in[0, \infty) \times [0, \infty) : L_{u}(p_{\lambda}, q_{\lambda}) = (0, 0),\lambda > 0\},$$
where $L_{u}$ is given by (\ref{eq:3.2}).\\
Since $p_{\lambda}u^{+}+q_{\lambda}u^{-} \in \mathcal{N}_{\lambda}$ , by assumption $(V_{2})$, (\ref{eq:3.7}) and (\ref{eq:3.8}), we have
\begin{align*}
  \displaystyle p_{\lambda}^{N}\|u^{+}\|^{N}+q_{\lambda}^{N}\|u^{-}\|^{N}
  &=\lambda\int_{B}f(x,p_{\lambda}u^{+}+q_{\lambda}u^{-})(p_{\lambda}u^{+}+q_{\lambda}u^{-})dx\\
  &\geq  \lambda \theta C_{5}p_{\lambda}^{\theta}\int_{B}|u^{+}|^{\theta}dx+\lambda \theta C_{5}q_{\lambda}^{\theta}\int_{B}|u^{-}|^{\theta}dx-\lambda \theta C_{6}|B|.
  \end{align*}
Since $\theta>N$, the set  $\mathcal{T}_{u}$ is bounded. Therefore, if $\{\lambda_{n}\}\subset (0,\infty)$
satisfies $\lambda_{n}\rightarrow \infty$ as $n\rightarrow\infty$, then up to subsequence, there exists $\bar{p}, \bar{q}>0$,
such that $p_{\lambda_{n}}\rightarrow \bar{p}$ and $q_{\lambda_{n}}\rightarrow \bar{q}$ .\\
We Claim that $\bar{p}=\bar{q}=0$. We proceed by contradiction and suppose that  $\bar{p}>0$ and $\bar{q}>0$.
For each $n\in \N$, $p_{\lambda_{n}}u^{+}+q_{\lambda_{n}}u^{-}\in\mathcal{N}_{\lambda_{n}}$. So,
\begin{align}\label{eq:3.17}
  \displaystyle \|p_{\lambda_{n}}u^{+}+q_{\lambda_{n}}u^{-}\|^{N}&=
 \lambda_{n}\int_{B}f(p_{\lambda_{n}}u^{+}+q_{\lambda_{n}}u^{-})(p_{\lambda_{n}}u^{+}+q_{\lambda_{n}}u^{-})dx.
\end{align}
It should be noted that $p_{\lambda_{n}}u^{+}\rightarrow \bar{p}u^{+}$ and $q_{\lambda_{n}}u^{-}\rightarrow \bar{q}u^{-}$ in $\mathbf{E}.$\\

On one hand, $\lambda_{n}\rightarrow 0$ as $n\rightarrow\infty$ and $\{p_{\lambda_{n}}u^{+}+q_{\lambda_{n}}u^{-}\}$ is bounded in $\mathbf{E}$. On the other hand, from
(\ref{eq:3.17}),we have
\begin{align*}
  \displaystyle \int_{B}|\nabla(\bar{p}u^{+}+\bar{q}u^{-})|^{N}dx  =
  \left(\lim_{n\rightarrow \infty}\lambda_{n}\right)\lim_{n\rightarrow \infty}\int_{B}f(p_{\lambda_{n}}u^{+}+q_{\lambda_{n}}u^{-})(p_{\lambda_{n}}u^{+}+q_{\lambda_{n}}u^{-})dx
\end{align*}
which is impossible.\\
Thus, $\bar{p}=\bar{q}=0$, so, $p_{\lambda_{n}}\rightarrow 0$ and $q_{\lambda_{n}}\rightarrow 0$ as   $n\rightarrow\infty$. Finally,
by $(V_{2})$ and (\ref{eq:3.17}), we have
$$\displaystyle 0\leq c_{\lambda}=\inf_{\mathcal{N}_{\lambda}} \mathcal{J}_{\lambda}(u)\leq  \mathcal{J}_{\lambda}(p_{\lambda_{n}}u^{+}+q_{\lambda_{n}}u^{-})\rightarrow 0. $$
Consequently, $c_{\lambda}\rightarrow 0$ as $\lambda\rightarrow \infty$.
\begin{lemma}\label{lem5}
If $ u_{0}\in \mathcal{N}_{\lambda} $ satisfies $\mathcal{J}_{\lambda}(u_{0})=c_{\lambda}$, then $\displaystyle\mathcal{J}'_{\lambda}(u_{0})=0.$
\end{lemma}
$\mathbf{Proof}$.We proceed by contradiction. We assume that $\displaystyle\mathcal{J}'_{\lambda}(u_{0})\neq 0$. By the continuity of $\mathcal{J}'_{\lambda}$,
there exists $\iota, \delta\geq 0$ such that
\begin{equation}\label{eq:3.18}
\displaystyle \|\mathcal{J}'_{\lambda}(v)\|_{\mathbf{E}^{\ast}}\geq\iota \mbox{ for all } \|v-u_{0}\|\leq 3\delta.
\end{equation}
Choose $\displaystyle \tau \in(0, \min \{\frac{1}{4}, \frac{\delta}{4\|u_{0}\|}\})$. Let $\displaystyle D=\left(1-\tau,1+\tau\right)\times \left(1-\tau,1+\tau\right)$ and define $g:D\rightarrow\mathbf{E}$, by

$$\displaystyle g(\rho,\vartheta)=\rho u_{0}^{+}+\vartheta u_{0}^{-},  (\rho,\vartheta) \in D$$
By virtue of $u_{0} \in \mathcal{N}_{\lambda}$, $\mathcal{J}_{\lambda}(u_{0})=c_{\lambda}$ and Lemma \ref{lem1}, it is easy to see that
\begin{equation}\label{eq:3.19}
\displaystyle\bar{c_{\lambda}}:=\max_{\partial D} \mathcal{J}_{\lambda}\circ g<c_{\lambda}.
\end{equation}
 Let $\epsilon:=\min\{\frac{c_{\lambda}-\bar{c_{\lambda}}}{3}, \frac{\iota\delta}{8}\}$, $S_{r}:=B(u_{0},r),r\geq0$
and $\displaystyle\mathcal{J}_{\lambda}^{a}:=\mathcal{J}_{\lambda}^{-1}(]-\infty,a]).$
 According to the Quantitative Deformation Lemma $[\cite{Wi}, \mbox{ Lemma }2.3]$,
there exists a deformation $\eta \in C\left([0,1]\times g(D), \mathbf{E}\right)$ such that:
 \begin{itemize}
   \item [$(1)$] $\eta(1, v)=v,$ if $v\not\in \mathcal{J}_{\lambda}^{-1}([ c_{\lambda}-2\epsilon,c_{\lambda}+2\epsilon])\cap S_{2\delta}$
   \item[$(2)$] $\eta\left(1, \mathcal{J_{\lambda}}^{c_{\lambda}+\epsilon }\cap S_{\delta}\right)\subset \mathcal{J_{\lambda}}^{c_{\lambda}-\epsilon}$,
    \item[$(3)$] $\mathcal{J_{\lambda}}(\eta(1, v))\leq \mathcal{J}_{\lambda}(v)$, for all $v\in  \mathbf{E}. $
 \end{itemize}
By lemma \ref{lem1} $(ii)$, we have $\mathcal{J}_{\lambda}( g(\rho,\vartheta))\leq c_{\lambda}$. In addition, we have,$$\|g(s,t)-u_{0}\|=\|(\rho-1)u^{+}_{0}+(\vartheta-1))u^{-}_{0}\|\leq |\rho-1|\|u^{+}_{0}\|+|\vartheta-1|\|u^{-}_{0}\|\leq2\tau \|u_{0}\|,$$ then $g(\rho,\vartheta)\in S_{\delta}$ for $(\rho,\vartheta)\in \bar{D}$.
Therefore, it follows from $(2)$ that
\begin{equation}\label{eq:3.20}
\max_{(\rho,\vartheta)\in \bar{D}}\mathcal{J}_{\lambda}(\eta(1, g(\rho,\vartheta)))\leq c_{\lambda}-\epsilon.
\end{equation}
In the following, we prove that $\eta(1, g(D))\cap\mathcal{N}_{\lambda}$ is nonempty. And in this case it contradicts (\ref{eq:3.20}) due to the definition
of $c_{\lambda}$.
 To do this, we first define
$$\bar{g}(\rho,\vartheta):=\eta(1, g(\rho,\vartheta)),$$
\begin{align*}
  \Upsilon_{0}(\rho,\vartheta)& = (\langle \mathcal{J}'_{\lambda}(g(\rho,\vartheta)), u_{0}^{+} \rangle, \langle \mathcal{J}'_{\lambda}(g(\rho,\vartheta)), u_{0}^{-} \rangle) \\
   &= (\langle \mathcal{J}'_{\lambda}(\rho u_{0}^{+}+\vartheta u_{0}^{-}), u_{0}^{+} \rangle, \langle \mathcal{J}'_{\lambda}(\rho u_{0}^{+}+\vartheta u_{0}^{-}), u_{0}^{-} \rangle) \\
  &:=(\varphi_{u_{0}}^{1}(\rho,\vartheta),\varphi_{u_{0}}^{2}(\rho,\vartheta))
\end{align*}
and
$$\Upsilon_{1}(\rho,\vartheta):=(\frac{1}{\rho}\langle \mathcal{J}'_{\lambda}(\bar{g}(\rho,\vartheta)),(\bar{g}(\rho,\vartheta))^{+} \rangle,
\frac{1}{\vartheta} \langle \mathcal{J}'_{\lambda}(\bar{g}(\rho,\vartheta)),(\bar{g}(\rho,\vartheta))^{-} \rangle).$$
Moreover, a simple calculation, shows that
\begin{align*}
  \displaystyle\frac{\varphi_{u_{0}}^{1}(\rho,\vartheta)}
  {\partial \rho}\bigg|_{(1,1)} & =(N-1)\|u_{0}^{+}\|^{N}
  -\lambda \int_{B}f'(x,u_{0}^{+})|u_{0}^{+}|^{2}dx\\
  &=(N-1)\lambda\int_{B}f(u_{0}^{+})u_{0}^{+}dx
  -\lambda \int_{B}f'(x,u_{0}^{+})|u_{0}^{+}|^{2}dx
\end{align*}
and $$  \displaystyle\frac{\varphi_{u_{0}}^{1}(\rho,\vartheta)}
  {\partial \vartheta}\bigg|_{(1,1)}=0.$$
  In the same manner, $$\displaystyle\frac{\varphi_{u_{0}}^{2}(\rho,\vartheta)}
  {\partial \rho}\bigg|_{(1,1)}=0$$
  and
  \begin{align*}
  \displaystyle\frac{\varphi_{u_{0}}^{2}(\rho,\vartheta)}
  {\partial \vartheta}\bigg|_{(1,1)}&=(N-1)\lambda\int_{B}f(x,u_{0}^{-})u_{0}^{-}dx
  -\lambda \int_{B}f'(x,u_{0}^{-})|u_{0}^{-}|^{2}dx
\end{align*}
Let $$\displaystyle \mathrm{J}=\left(
          \begin{array}{cc}
            \frac{\varphi_{u_{0}}^{1}(\rho,\vartheta)}
  {\partial \rho}\bigg|_{(1,1)} & \frac{\varphi_{u_{0}}^{2}(\rho,\vartheta)}
  {\partial \rho}\bigg|_{(1,1)}\\
           \frac{\varphi_{u_{0}}^{1}(\rho,\vartheta)}
  {\partial \vartheta}\bigg|_{(1,1)} & \frac{\varphi_{u_{0}}^{2}(\rho,\vartheta)}
  {\partial \vartheta}\bigg|_{(1,1)} \\
          \end{array}
        \right).
$$
Then we have $\det \mathrm{J}\neq 0$. Therefore, the point $(0,1)$ is the unique isolated zero of the $C^{1}$ function $\Upsilon_{0}$. By using the Brouwer's degree in $\R^{2}$, we deduce that $\mbox{deg}(\Upsilon_{0}, D, 0)=1.$\\
Now, it follows from (\ref{eq:3.20}) and $(1)$ that $g(\rho,\vartheta)=\bar{g}(\rho,\vartheta)$ on $\partial D$.
For the boundary dependence of Brouwer's degree ( see [\cite{Du}, Theorem 4.5]), there holds $\mbox{deg}(\Upsilon_{1}, D, 0)=\mbox{deg}(\Upsilon_{0}, D, 0)=1.$
Therefore, there exists some $(\overline{\rho},\overline{\vartheta})\in D$ such that
$$\eta(1, g(\overline{\rho},\overline{\vartheta})\in \mathcal{N}_{\lambda}.$$
This finish the proof of the Lemma.
\begin{lemma}\label{lem6}
If $\upsilon$ is a least energy sign-changing solution of problem $(P_{\lambda})$, then $\upsilon$ has exactly two nodal domains
\end{lemma}
\textbf{Proof. }Assume by contradiction that $\upsilon=\upsilon_{1}+\upsilon_{2}+\upsilon_{3}$ satisfies
\begin{align*}
 \displaystyle\upsilon_{i}\neq 0, i& =1,2,3, \upsilon_{1}\geq0, \upsilon_{2}\leq0 , \mbox{ a.e. in } B \\
  B_{1}\cap B_{2}&=\emptyset, B_{1}:=\{ x \in B:\upsilon_{1}(x)>0\}, B_{2}:=\{ x \in B:\upsilon_{2}(x)<0\}
\end{align*}
$$\upsilon_{1}\bigg|_{B\setminus B_{1}\cup B_{2}}=\upsilon_{2}\bigg|_{B\setminus B_{2}\cup B_{1}}=\upsilon_{3}\bigg|_{B_{1}\cup B_{2}}=0,$$
and
\begin{equation}\label{eq:3.21}
 \langle \mathcal{J}'_{\lambda}(\upsilon), \upsilon_{i}\rangle=0 \mbox{ for }  i=1,2,3,
\end{equation}
Let $\nu=\upsilon_{1}+\upsilon_{2}$ and it is easy to see that $\nu^{+}=\upsilon_{1}$, $\nu^{-}=\upsilon_{2}$ and $\nu^{\pm}\neq0$. From Lemma (\ref{lem1}),
it follows that there exists a unique couple $(p_{\nu},q_{\nu})\in [ 0, \infty) \times[ 0, \infty)$ such that
$p_{\nu}\upsilon_{1}+q_{\nu}\upsilon_{2}\in \mathcal{N}_{\lambda}$. So,
$\mathcal{J_{\lambda}}(p_{\nu}\upsilon_{1}+q_{\nu}\upsilon_{2})\geq c_{\lambda}$. Moreover, using (\ref{eq:3.21}), we obtain that
$ \langle \mathcal{J}'_{\lambda}(\nu), \nu^{\pm}\rangle=0.$ Then, by Lemma (\ref{lem2}), we have
$0<p_{\nu}, q_{\nu}\leq 1$.\\

Now, combining (\ref{eq:3.21}), $(V_{3})$ and $(V_{4})$, we have that
\begin{align*}
  \displaystyle0& =\frac{1}{\theta}\langle \mathcal{J}'_{\lambda}(\upsilon), \upsilon_{3}\rangle =\frac{1}{\theta}\langle \mathcal{J}'_{\lambda}(\upsilon_{3}), \upsilon_{3}\rangle \\
  <& \mathcal{J}_{\lambda}(\upsilon_{3}),
\end{align*}
and
\begin{align*}
 \displaystyle c_{\lambda}\leq& \mathcal{J}_{\lambda}(p_{\nu}\upsilon_{1}+q_{\nu}\upsilon_{2})  \\
 =&\mathcal{J}_{\lambda}(p_{\nu}\upsilon_{1}+q_{\nu}\upsilon_{2})-\frac{1}{\theta}\langle \mathcal{J}'_{\lambda}p_{\nu}\upsilon_{1}+q_{\nu}\upsilon_{2}), p_{\nu}\upsilon_{1}+q_{\nu}\upsilon_{2}\rangle  \\
 &=(\frac{1}{N}-\frac{1}{\theta})p_{\nu}^{N}\|\upsilon_{1}\|^{N}
+(\frac{1}{N}-\frac{1}{\theta})q_{\nu}^{N}\|\upsilon_{2}\|^{N}
\\
 &+\lambda\int_{B}[ \frac{1}{\theta}f(x,p_{\nu}\upsilon_{1})(p_{\nu}\upsilon_{1})-F(x,p_{\nu}\upsilon_{2})]dx+\lambda\int_{B}[ \frac{1}{\theta}f(x,q_{\nu}\upsilon_{1})(p_{\nu}\upsilon_{2})-F(x,q_{\nu}\upsilon_{2})]dx\\
  &\leq \mathcal{J}_{\lambda}(\upsilon_{1}+\upsilon_{2})-\frac{1}{\theta}\langle \mathcal{J}'_{\lambda}\upsilon_{1}+\upsilon_{2}), \upsilon_{1}+\upsilon_{2}\rangle \\
& =\mathcal{J}_{\lambda}(\upsilon_{1}+\upsilon_{2})+\frac{1}{\theta}\langle \mathcal{J}'_{\lambda}(\upsilon),\upsilon_{3}\rangle\\
  &<\mathcal{J}_{\lambda}(\upsilon_{1}+\upsilon_{2})+\mathcal{J}_{\lambda}(\upsilon_{3})=\mathcal{J}_{\lambda}(\upsilon)= c_{\lambda},
 \end{align*}
which is a contradiction. Therefore, $\upsilon_{3}=0$ and $\upsilon$ has exactly two nodal domains.
\section{The subcritical case}
\begin{lemma}\label{lem7}
If $ \{u_{n}\}\subset \mathcal{N}_{\lambda} $ is a minimizing sequence for $c_{\lambda} $, then
 there exists some $u\in \mathbf{E}$ such that
$$\displaystyle\int_{B}f(u_{n}^{\pm})u_{n}^{\pm} dx\rightarrow \int_{B}f(u^{\pm})u^{\pm} dx $$
and
$$\int_{B}F(u_{n}^{\pm}) dx\rightarrow \int_{B}F(u^{\pm}) dx $$
\end{lemma}
\textbf{Proof.} we will only prove the first  result. Since the second limit is a direct consequence of the first one, we omit
it here.\\
Let sequence $\displaystyle\{u_{n}\}\subset\mathcal{N}_{\lambda} $ satisfy $\displaystyle\lim_{n\rightarrow \infty}\mathcal{J_{\lambda}}(u_{n}) =c_{\lambda}$.
It is clearly that $\{u_{n}\}$ is bounded by Lemma (\ref{lem3}). Then , up to a subsequence, there exists $u\in \mathbf{E}$ such that
\begin{align}
 \displaystyle u_{n}&\rightharpoonup u \mbox{ in } \mathbf{E},\nonumber \\
  u_{n} &\rightarrow u   \mbox{ in } L^{t} (B) \mbox{ for } t \in [1, \infty),\label{eq:4.1}\\
 u_{n} &\rightarrow u \mbox{ a.e. in  } B. \nonumber \\
 u_{n}^{\pm} &\rightharpoonup u^{\pm}   \mbox{ in } \mathbf{E},\nonumber \\
  u_{n}^{\pm} & \rightarrow u^{\pm} \mbox{ in }  L^{t} (B) \mbox{ for } t \in [1, \infty), \label{eq:4.2}\\
  u_{n}^{\pm} & \rightarrow u^{\pm}\mbox{ a.e. in  } B\nonumber
\end{align}
Note that by (\ref{eq:3.3}), we have
\begin{equation}\label{eq:4.3}
\displaystyle  f(x,u_{n}^{\pm}(x))u_{n}^{\pm}(x)\leq \epsilon |u_{n}^{\pm}(x)|^{N}+C_{2} |u_{n}^{\pm}(x)|^{s}
\exp (\alpha |u_{n}^{\pm}(x)|^{\gamma}):=h(u_{n}^{\pm}(x)),
\end{equation}
for all $\alpha>\alpha_{0}$ and $q>N.$ It is sufficient to prove that sequence $\{h(u_{n}^{\pm})\}$
is convergent in $L^{1}(B)$.\\
Choosing $a,a'>1$ with $\frac{1}{a}+\frac{1}{a'}=1$, we get that
\begin{equation}\label{eq:4.4}
\displaystyle |u_{n}^{\pm}|^{s}\rightarrow |u^{\pm}|^{s} \mbox{ in }  L^{a'} (B)
\end{equation}
Moreover, choosing $\alpha>0$ small enough  such that $\alpha a \left(\displaystyle\max_{n}\|u_{n}^{\pm}\|^{\gamma}\right)\leq \alpha_{N,\beta}$,
 we conclude from (\ref{eq:1.4}) that
\begin{equation}\label{eq:4.5}
  \displaystyle\int_{B}\exp\left(\alpha |u_{n}^{\pm}(x)|^{\gamma}\right)dx<\infty.
\end{equation}
Since $\exp\left(\alpha |u_{n}^{\pm}(x)|^{\gamma}\right)dx \rightarrow \exp\left(\alpha |u^{\pm}(x)|^{\gamma}\right)dx $, a.e. in $B$.
 From (\ref{eq:4.5}) and [\cite{Kav}, Lemma 4.8, chapter 1], we obtain that
\begin{equation}\label{eq:4.6}
 \displaystyle\exp \left(\alpha |u_{n}^{\pm}|^{\gamma}\right)dx\rightharpoonup \exp \left(\alpha |u^{\pm}|^{\gamma}\right)dx \mbox{ in }  L^{a}(B).
\end{equation}

Hence, by (\ref{eq:4.4}), (\ref{eq:4.6}) and [\cite{Kav}, Lemma 4.8, chapter 1] again, we conclude that
$$\displaystyle\int_{B}f(u_{n}^{\pm})u_{n}^{\pm} dx\rightarrow \int_{B}f(u^{\pm})u^{\pm} dx $$

\begin{lemma}\label{lem8}
There exists some $ \upsilon \in \mathcal{N}_{\lambda} $ such that $\mathcal{J}_{\lambda}(\upsilon)=c_{\lambda}$.
\end{lemma}
\textbf{Proof.} Let $\{\upsilon_{n}\}\subset \mathcal{N}_{\lambda} $ be a sequence such that $\displaystyle\lim_{n\rightarrow\infty}\mathcal{J_{\lambda}}(\upsilon_{n})=c_{\lambda}.$
It is clearly that $\{\upsilon_{n}\}$ is bounded by Lemma (\ref{lem3}). Then, up to a subsequence, there exists $\upsilon \in \mathbf{E}$ such that
\begin{align}
   \displaystyle\upsilon_{n}^{\pm} &\rightharpoonup \upsilon^{\pm}   \mbox{ in } \mathbf{E},\nonumber \\
  \upsilon_{n}^{\pm} & \rightarrow \upsilon^{\pm}\mbox{ in }  L^{t} (B) \mbox{ for } t \in [1, \infty), \label{eq:4.7}\\
  \upsilon_{n}^{\pm} & \rightarrow \upsilon^{\pm}\mbox{ a.e. in  } B\nonumber
\end{align}
We claim that $\upsilon^{+}\neq 0$ and $\upsilon^{-}\neq 0$. Suppose, by contradiction, $\upsilon^{+}= 0$ .
From the definition of $\mathcal{N}_{\lambda} $, (\ref{eq:4.7}), (\ref{eq:4.3}) and Lemma (\ref{lem7}), we have that
$\displaystyle\lim_{n\rightarrow\infty}\| \upsilon_{n}^{+}\|^{N}=0$, which contradicts Lemma (\ref{lem3}). Hence, $\upsilon^{+}\neq 0$
and $\upsilon^{-}\neq 0$.\\
From the lower semi continuity of norm and (\ref{eq:4.7}), it follows that
\begin{equation}\label{eq:4.8}
  \displaystyle\langle\mathcal{J}'_{\lambda}(\upsilon), \upsilon^{\pm}\rangle \leq\lim_{n\rightarrow\infty}\langle\mathcal{J}'_{\lambda}(\upsilon_{n}), \upsilon_{n}^{\pm}\rangle =0.
\end{equation}
Then, Lemma (\ref{lem2}) implies that there exists $(p_{\upsilon},q_{\upsilon} )\in (0,1]\times(0,1]$ such that
$p_{\upsilon}\upsilon^{+}+q_{\upsilon}\upsilon^{-}\in \mathcal{N}_{\lambda}$. Thus, by $(V_{2})$, $\lambda\geq0$ and Lemma (\ref{lem7}), we get that
\begin{align}\label{eq:4.9}
 \displaystyle c_{\lambda} \leq\mathcal{J_{\lambda}}(p_{\upsilon}\upsilon^{+}+q_{\upsilon}\upsilon^{-})=&\mathcal{J_{\lambda}}(p_{\upsilon}\upsilon^{+}+q_{\upsilon}\upsilon^{-})
-\frac{1}{\theta}\langle\mathcal{J}'_{\lambda}(p_{\upsilon}\upsilon^{+}+q_{\upsilon}\upsilon^{-}), p_{\upsilon}\upsilon^{+}+q_{\upsilon}\upsilon^{-}\rangle \nonumber \\
  \leq&\mathcal{J}(\upsilon)
-\frac{1}{\theta}\langle\mathcal{J}'_{\lambda}(\upsilon), \upsilon\rangle   \\
   \leq& \lim_{n\rightarrow \infty}\left[\mathcal{J_{\lambda}}(\upsilon_{n})
-\frac{1}{\theta}\langle\mathcal{J}'_{\lambda}(\upsilon_{n}), \upsilon_{n}\rangle \right] \nonumber\\
   =& \lim_{n\rightarrow \infty}\mathcal{J_{\lambda}}(\upsilon_{n})=c_{\lambda}\cdot\nonumber
\end{align}
Noticing that if $p_{\upsilon}<1$ or $q_{\upsilon}<1$, then the inequality (\ref{eq:4.9}) is strict. Hence, by bringing together (\ref{eq:4.8}) and (\ref{eq:4.9}),
we conclude that $p_{\upsilon}=q_{\upsilon}=1$ and $\upsilon\in \mathcal{N}_{\lambda}$ satisfying $\mathcal{J}(\upsilon)=c_{\lambda}$.\\
\textbf{Proof of Theorem \ref{th1.2}. }From Lemma \ref{lem5}, Lemma \ref{lem6} and Lemma \ref{lem8}, we deduce that
$\upsilon$ is a least energy sign-changing solution form problem $(P_{\lambda})$ with exactly tow nodal domains.

\section{The critical case}

\begin{lemma}\label{lem10}There exists $\lambda^{\ast}>0$ such that
if $\lambda\geq\lambda^{\ast}$, and $\{\upsilon_{n}\}\subset \mathcal{N}_{\lambda}$
is a minimizing sequence for $c_{\lambda}$,   then there exists some $ \upsilon \in \mathcal{N}_{\lambda} $ such that $\mathcal{J}_{\lambda}(\upsilon)=c_{\lambda}$.
\end{lemma}
\textbf{Proof.} Let $\{\upsilon_{n}\}\subset \mathcal{N}_{\lambda} $ be a sequence such that $\displaystyle\lim_{n\rightarrow\infty}\mathcal{J_{\lambda}}(v_{n})=c_{\lambda}.$ We have \begin{align*}\label{eq:5.5} \mathcal{J}_{\lambda}(\upsilon_{n})\rightarrow c_{\lambda}~~\mbox{and}~~\langle\mathcal{J}'_{\lambda}(\upsilon_{n}),\varphi \rangle\rightarrow 0, \forall \varphi\in\mathbf{E}\end{align*}that is
\begin{equation}\label{eq:5.1}
\mathcal{J_{\lambda}}(\upsilon_{n})=\frac{1}{N}\|\upsilon_{n}\|^{N}-\int_{B}F(x,\upsilon_{n})dx \rightarrow c_{\lambda} ,~~n\rightarrow +\infty
\end{equation}and
\begin{equation}\label{eq:5.2}
| \langle\mathcal{J}'_{\lambda}(\upsilon_{n}),\varphi \rangle|=\Big|\int_{B}\omega(x)|\nabla \upsilon_{n}|^{N-2}\nabla \upsilon_{n}.\nabla \varphi dx -\int_{B}f(x,\upsilon_{n})\varphi  dx
\Big|\leq \varepsilon_{n}\|\varphi\|,
\end{equation}
for all $\varphi \in \mathbf{E}$, where $\varepsilon_{n}\rightarrow0$, as $n\rightarrow +\infty$.\\
By lemma \ref{lem3}, $\upsilon_{n}$ is bounded in $\mathbf{E}$.
 Furthermore, we have from (\ref{eq:5.2}) and $(V_{2})$, that
\begin{equation}\label{eq:5.3}
0<\int_{B} f(x,u_{n})u_{n}\leq C
 \end{equation}
and
 \begin{equation*}\label{eq:5.6}
0<\int_{B} F(x,u_{n})\leq C.
 \end{equation*}
Since by Lemma \ref{lem2}, we have
\begin{equation}\label{eq:5.4}
f(x,u_{n})\rightarrow f(x,u) ~~\mbox{in}~~L^{1}(B) ~~as~~ n\rightarrow +\infty,
 \end{equation}
then, it follows from $(H_{2})$ and the generalized Lebesgue dominated convergence Theorem that
\begin{equation}\label{eq:5.5}
F(x,u_{n})\rightarrow F(x,u) ~~\mbox{in}~~L^{1}(B) ~~as~~ n\rightarrow +\infty.
 \end{equation}
Arguing as Lemma \ref{lem8}, we have that, up to a subsequence,
\begin{align}
 \displaystyle \upsilon_{n}&\rightharpoonup u \mbox{ in } \mathbf{E},\nonumber \\
  \upsilon_{n} &\rightarrow u   \mbox{ in } L^{t} (B) \mbox{ for } t \in [1, \infty),\label{eq:5.6}\\
 \upsilon_{n} &\rightarrow u \mbox{ a.e. in  } B. \nonumber \\
 \upsilon_{n}^{\pm} &\rightharpoonup u^{\pm}   \mbox{ in } \mathbf{E},\nonumber \\
  \upsilon_{n}^{\pm} & \rightarrow u^{\pm} \mbox{ in }  L^{t} (B) \mbox{ for } t \in [1, \infty), \label{eq:5.7}\\
  \upsilon_{n}^{\pm} & \rightarrow u^{\pm}\mbox{ a.e. in  } B\nonumber
\end{align}

for some $\upsilon \in \mathbf{E}$.\\Noticing that, according to lemma \ref{lem4}, there exists $\lambda^{*} >0$ such that for all $\lambda> \lambda^{*}$, we get $$c_{\lambda}<\displaystyle\frac{1}{N}(\frac{\alpha_{N,\beta}}{\alpha_{0}})^{\frac{N}{\gamma}}\cdot$$ In the sequel, the results that are valid for $\upsilon_{n}~~\mbox{and}~~\upsilon$, are also valid for $\upsilon^{\pm}_{n}~~\mbox{and}~~\upsilon^{\pm}$. Next, we are going to make some Claims.
\medskip\\
\noindent {\bf\it\underline{Claim 1}.}  $\nabla \upsilon_{n}(x)\rightarrow\nabla \upsilon(x) ~~a.e.$ in $ B$ and $\upsilon$ is a solution of the problem $(P_{\lambda})$.\\ Indeed, for any $\xi>0$, let $\mathcal{A_{\eta}}=\{x\in B, |\upsilon_{n}-\upsilon|\geq \xi\}$. For all $t\in\mathbb{R}$, for all positive $c>0$, we have $$ct\leq e^{t}+c^{2}.$$ It follows that for $t=\alpha_{N,\beta}\big(\frac{|\upsilon_{n}-\upsilon|}{\|\upsilon_{n}-\upsilon\|}\big)^{\gamma}$, $c=\frac{1}{\alpha_{N,\beta}}\|\upsilon_{n}-\upsilon\|^{\gamma}$, we get
     $$\begin{array}{rclll}
|\upsilon_{n}-\upsilon|^{\gamma}&\leq&{\displaystyle e^{\alpha_{N,\beta}\big(\frac{|\upsilon_{n}-\upsilon|}{\|\upsilon_{n}-\upsilon\|}\bigl)^{\gamma}}+\frac{1}{\alpha_{N,\beta}^{2}}\|\upsilon_{n}-\upsilon\|^{2\gamma}}
\\
 &\leq&{\displaystyle e^{\alpha_{N,\beta}\bigl(\frac{|\upsilon_{n}-\upsilon|}{\|\upsilon_{n}-\upsilon\|}\bigl)^{\gamma}}+C_{1}(N)},
    \end{array}$$
    where $C_{1}(N)$ is a constant depending only on $N$ and the upper bound of $\|\upsilon_{n}\|$. So, if we denote by $\mathcal{L}(\mathcal{A_{\xi}})$ the Lebesgue measure of the set $\mathcal{A_{\xi}}$, we obtain
      $$\begin{array}{rclll}
\displaystyle\mathcal{L}(\mathcal{A_{\xi}})=\displaystyle\int_{\mathcal{A_{\xi}}}e^{|\upsilon_{n}-\upsilon|^{\gamma}}e^{-|\upsilon_{n}-\upsilon|^{\gamma}}dx&\leq&
e^{-\xi^{\gamma}}\displaystyle\int_{\mathcal{A_{\xi}}}\exp\Big( e^{\alpha_{N,\beta}\bigl({\normalsize \frac{|\upsilon_{n}-\upsilon|}{\|\upsilon_{n}-\upsilon\|}}\bigl)^{\gamma}}+C_{1}(N)\Big)dx
\\
 &\leq&\displaystyle e^{-\xi^{\gamma}}e^{C_{1}(N)}\displaystyle\int_{B} \exp\Big(\alpha_{N,\beta} \bigl({\normalsize \frac{|\upsilon_{n}-\upsilon|}{\|\upsilon_{n}-\upsilon\|}}\bigl)^{\gamma}\Big)dx\\
 &\leq& e^{-\xi^{\gamma}}C_{2}(N)\rightarrow0 ~~\mbox{as}~~\xi\rightarrow+\infty,
    \end{array}$$
  where $C_{2}(N)$ is a positive constant depending only on $N$ and the upper bound of  $\|\upsilon_{n}\|$. It follows that
    \begin{equation}\label{eq:5.8}
    \int_{\mathcal{A_{\xi}}}|\nabla \upsilon_{n}-\nabla \upsilon|dx\leq Ce^{-\frac{1}{2}\xi^{\gamma}}\Big(\int_{B}|\nabla \upsilon_{n}-\nabla \upsilon|^{2}\omega (x)dx\Big)^{\frac{1}{2}}\rightarrow0 ~~\mbox{as}~~\xi\rightarrow+\infty\cdot
    \end{equation}
   We define for  $\xi>0$, the following truncation function
$$T_{\xi}(s):=\displaystyle \left\{
      \begin{array}{rcll}
&s& \mbox{ if }|s| <\xi\\
       &\xi\frac{s}{|s|}& \mbox{ if } |s|\geq\xi.\\
 \end{array}
 \right.$$
If we take $\varphi=T_{\xi}(\upsilon_{n}-\upsilon)\in \mathbf{E}$, in (\ref{eq:5.2}) then with $\nabla \varphi=\chi_{_{\mathcal{A_{\xi}}}}\nabla(\upsilon_{n}-\upsilon)$,  we obtain
$$\begin{array}{rclll}
\displaystyle\Big|\int_{B\setminus \mathcal{A_{\xi}}}\omega(x)\big(|\nabla \upsilon_{n}|^{N-2}\nabla
\upsilon_{n}-|\nabla \upsilon|^{N-2}\nabla \upsilon\big).(\nabla \upsilon_{n}-\nabla
\upsilon)dx\Big| &\leq&\\
\displaystyle\Big|\int_{B\setminus \mathcal{A_{\xi}}}\omega(x)|\nabla \upsilon|^{N-2}\nabla \upsilon.(\nabla \upsilon_{n}-\nabla \upsilon)dx\Big|+
\int_{B}f(x,\upsilon_{n})T_{\xi}(\upsilon_{n}-\upsilon)dx +\varepsilon_{n}\|\upsilon_{n}-\upsilon\|&\leq&\\
\displaystyle\Big|\int_{B}\omega(x)|\nabla \upsilon|^{N-2}\nabla \upsilon.(\nabla \upsilon_{n}-\nabla \upsilon)dx\Big|+
\int_{B}f(x,\upsilon_{n})T_{\xi}(\upsilon_{n}-\upsilon)dx +\varepsilon_{n}\|\upsilon_{n}-\upsilon\|
\end{array}$$ where $\varepsilon_{n}\rightarrow0~~\mbox{as}~~n\rightarrow+\infty$.\\
Since $\upsilon_{n}\rightharpoonup \upsilon ~~\mbox{weakly}$, then $\displaystyle\int_{B}\omega(x)|\nabla \upsilon|^{N-2}\nabla \upsilon.(\nabla \upsilon_{n}-\nabla \upsilon)\rightarrow0$. Moreover, by (\ref{eq:5.4}) and   the Lebesgue dominated convergence Theorem, we get
$$ \displaystyle\int_{B}f(x,\upsilon_{n})T_{\xi}(\upsilon_{n}-\upsilon)dx \rightarrow0~~\mbox{as}~~n\rightarrow+\infty.$$
Using the well known inequality,
$$\langle|x|^{N-2}x-|y|^{N-2}y,x-y\rangle~ \geq 2^{2-N}|x-y|^{N}~\forall ~~x,y \in \mathbb{R}^{N},~~N\geq 2,$$ $\langle \cdot, \cdot\rangle$ is the inner product in
$\mathbb{R^{N}}$, one has
$$\int_{B\setminus\mathcal{A_{\xi}}}\omega(x)|\nabla \upsilon_{n}-\nabla \upsilon|^{N}dx\rightarrow0\cdot$$
Therefore,
\begin{equation}\label{eq:5.9}
\int_{B\setminus\mathcal{A_{\xi}}}|\nabla \upsilon_{n}-\nabla \upsilon|dx\leq \Big(\int_{B\setminus\mathcal{A_{\xi}}}\omega(x)|\nabla \upsilon_{n}-\nabla \upsilon|^{N}dx\Big)^{\frac{1}{N}}\Big(\mathcal{L}(B\setminus\mathcal{A_{\xi}})\Big)^{\frac{1}{N'}}\rightarrow0~~\mbox{as}~~n\rightarrow+\infty.
\end{equation}
From (\ref{eq:5.8}) and (\ref{eq:5.9}), we deduce that $$\int_{B}|\nabla \upsilon_{n}-\nabla \upsilon|dx\rightarrow0~~\mbox{as}~~n\rightarrow+\infty.$$
Therefore,  $\nabla \upsilon_{n}(x)\rightarrow\nabla u(x) ~~a.e.$ in $B$.

On the other hand,$$\Big(|\nabla \upsilon_{n}|^{N-2}\nabla \upsilon_{n}\big)~~\mbox{is bounded in}~~(L^{\frac{N}{N-1}}(B,\omega))^{N}.$$ Then, up to subsequence, we can assume that\begin{equation}\label{eq:5.10}|\nabla \upsilon_{n}|^{N-2}\nabla \upsilon_{n}\rightharpoonup |\nabla
\upsilon|^{N-2}\nabla \upsilon ~~\mbox{weakly in
}~~(L^{\frac{N}{N-1}}(B,\omega))^{N}\cdot\end{equation}
Therefore, passing to the limit in (\ref{eq:5.2}) and using (\ref{eq:5.4}), (\ref{eq:5.10}), the convergence everywhere of the gradient, we obtain that $\upsilon$ is a solution of problem $(P_{\lambda})$. Claim 1 is proved.\\
\noindent {\bf\it\underline{Claim 2}.}  $\upsilon^{+}\neq 0$ and $\upsilon^{-}\neq 0$. Suppose, by contradiction, $\upsilon^{+}= 0$ . Therefore,
 $\displaystyle\int_{B} F(x,\upsilon_{n})dx~~\rightarrow 0$ and consequently we get\begin{equation}\label{eq:5.11}
\frac{1}{N}\|\upsilon_{n}\|^{N}\rightarrow c_{\lambda}<\frac{1}{N}(\frac{\alpha_{N,\beta}}{\alpha_{0}})^{\frac{N}{\gamma}}\cdot\end{equation}
 First, we claim that there exists $q>1$ such that
\begin{equation}\label{eq:5.12}
\int_{B}|f(x,\upsilon_{n})|^{q}dx\leq C.
\end{equation}
 By (\ref{eq:5.2}), we have
$$\Big| \|\upsilon_{n}\|^{N}-\int_{B}f(x,\upsilon_{n})\upsilon_{n}dx\Big|\leq C\varepsilon_{n}.$$
 So
$$\|\upsilon_{n}\|^{N}\leq C\varepsilon_{n}+\Big(\int_{B}|f(x,\upsilon_{n})|^{q}dx\Big)^{\frac{1}{q}}\Big(\int_{B}|\upsilon_{n}|^{q'}dx\Big)^{\frac{1}{q'}},$$
where $q'$ is the conjugate of $q$. Since $(\upsilon_{n})$ converge to 0 in $L^{q'}(B)$
$$\displaystyle\lim_{n\rightarrow +\infty}\|\upsilon_{n}\|^{N}=0.$$
According to Lemma \ref{lem3}, this result cannot occur. Now
for the proof of the claim (\ref{eq:5.12}), since $f$ has critical growth, for every
$\varepsilon>0$ and $q>1$ there exists $t_{\varepsilon}>0$ and $C>0$
such that for all $|t|\geq
t_{\varepsilon}$, we have
\begin{equation*}\label{eq:5.13}
|f(x,t)|^{q}\leq
Ce^{\alpha_{0}(\varepsilon+1) t^{\gamma}}.
\end{equation*}
Consequently,
$$\begin{array}{rcll}
\displaystyle\int_{B}|f(x,\upsilon_{n})|^{q}dx&=&\displaystyle\int_{\{|\upsilon_{n}|\leq
t_{\varepsilon}\}}|f(x,\upsilon_{n})|^{q}dx
+\int_{\{|\upsilon_{n}|>t_{\varepsilon}\}}|f(x,\upsilon_{n})|^{q}dx\\
 &\leq &\displaystyle\omega_{N-1} \max_{B\times [-t_{\varepsilon},t_{\varepsilon}]}|f(x,t)|^{q}+ C\int_{B}e^{\alpha_{0}(\varepsilon+1)|\upsilon_{n}|^{\gamma}}dx.
 \end{array}$$
Since $Nc_{\lambda}<\displaystyle(\frac{\alpha_{N,\beta}}{\alpha_{0}})^{\frac{N}{\gamma}}$, there exists
$\eta\in(0,\frac{1}{N})$ such that
$Nc_{\lambda}=\displaystyle(1-\eta)\displaystyle(\frac{\alpha_{N,\beta}}{\alpha_{0}})^{\frac{N}{\gamma}}$.
On the other hand, $\|\upsilon_{n}\|^{\gamma}\rightarrow
(Nc_{\lambda})^{\frac{\gamma}{N}}$, so there exists $n_{\eta}>0$ such that for
all $n\geq n_{\eta}$, we get $\|\upsilon_{n}\|^{\gamma}\leq
(1-\eta)\frac{\alpha_{N,\beta}}{\alpha_{0}}$. Therefore,
$$\alpha_{0}(1+\varepsilon)(\frac{|\upsilon_{n}|}{\|\upsilon_{n}\|})^{\gamma}\|\upsilon_{n}\|^{\gamma}\leq
(1+\varepsilon)(1-\eta)\alpha_{N,\beta}\cdot$$
 We choose $\varepsilon >0$ small enough to get
 $$\alpha_{0}(1+\varepsilon) \|\upsilon_{n}\|^{\gamma}\leq \alpha_{N,\beta}\cdot $$
So, the second integral is
uniformly bounded in view of (\ref{eq:1.4}) and the claim is
proved.

\medskip

Since ($\upsilon_{n}$) is bounded, up to a subsequence, we can assume that $\|\upsilon_{n}\|\rightarrow \rho>0.$ We affirm that $\mathcal{J_{\lambda}}(\upsilon)=c_{\lambda}$. Indeed, by $(V_{2})$ and claim 2, we have
\begin{equation}\label{eq:5.13}
\mathcal{J_{\lambda}}(\upsilon)= \frac{1}{N}\int_{B}[f(x,\upsilon)\upsilon-NF(x,\upsilon)]dx\geq0.
\end{equation}
Now, using the lower semi continuity of the norm and (\ref{eq:5.5}) , we get,
 $$\mathcal{J_{\lambda}}(\upsilon)\leq \frac{1}{N}\liminf_{n\rightarrow\rightarrow\infty} \|\upsilon_{n}\|^{N}-\int_{B}F(x,\upsilon)dx=c_{\lambda}.$$
Suppose that
\begin{equation*}\label{eq:5.16}\mathcal{J_{\lambda}}(\upsilon)<c_{\lambda}.\end{equation*} Then
\begin{equation}\label{eq:5.14}\|\upsilon\|^{N}<\rho^{N}.\end{equation} In addition,
\begin{equation}\label{eq:5.15}
\displaystyle\frac{1}{N}\lim_{n\rightarrow+\infty} \|\upsilon_{n}\|^{N}=\big(c_{\lambda}+\int_{B}F(x,\upsilon)dx\big),
\end{equation}
which means that $$\rho^{N}=N\Big (c_{\lambda}+\int_{B}F(x,\upsilon)dx\Big).$$
Set
$$u_{n}=\frac{\upsilon_{n}}{\|\upsilon_{n}\|}\quad \mbox{and}\quad u=\frac{\upsilon}{\rho}\cdot$$
We have $\|u_{n}\|=1$, $u_{n}\rightharpoonup u$  in $\mathbf{E}$, $u\not\equiv 0$ and $\|u\|<1$. So, by Lemma \ref{Lionstype}, we get
$$\displaystyle\sup_{n}\int_{B}e^{p~\alpha_{N,\beta}
|u_{n}|^{\gamma}}dx< +\infty,~$$  provided
$1<p<\bigl(1-\|u\|^{N}\bigl)^{-\frac{\gamma}{N}}$.\\

\noindent By (\ref{eq:5.5}) and (\ref{eq:5.15}), we have the following equality
$$Nc_{\lambda}-N\mathcal{J_{\lambda}}(\upsilon)=\rho^{N}-\|\upsilon\|^{N}.$$ From (\ref{eq:5.13}), Lemma \ref{lem7} and  the last equality, we obtain\begin{equation}\label{eq:5.16}
\rho^{N}\leq Nc_{\lambda}+\|\upsilon\|^{N}<(\frac{\alpha_{N,\beta}}{\alpha_{0}})^{\frac{N}{\gamma}}+\|\upsilon\|^{N}.\end{equation}

\noindent Since $$\rho^{\gamma}=\Big(\frac{\rho^{N}-\|\upsilon\|^{N}}{1-\|u\|^{N}}\Big)^{\frac{1}{(N-1)(1-\beta)}},$$ we deduce from (\ref{eq:5.16}) that
\begin{equation}\label{eq:5.17}
\rho^{\gamma}<\Big(\frac{(\frac{\alpha_{N,\beta}}{\alpha_{0}})^{\frac{N}{\gamma}}}{1-\|u\|^{N}}\Big)^{\frac{1}{(N-1)(1-\beta)}}.\end{equation}

On one hand, we have this estimate $\displaystyle\int_{B}|f(x,\upsilon_{n})|^{q}dx<C$. Indeed, since $f$ has critical growth, for every
$\varepsilon>0$ and $q>1$ there exists $t_{\varepsilon}>0$ and $C>0$
such that for all $|t|\geq
t_{\varepsilon}$, we have
\begin{equation*}
|f(x,t)|^{q}\leq
Ce^{\alpha_{0}(\varepsilon+1) t^{\gamma}}.
\end{equation*}
So,
$$\begin{array}{rclll}
\displaystyle\int_{B}|f(x,\upsilon_{n})|^{q}dx&=&\displaystyle\int_{\{|\upsilon_{n}|\leq
t_{\varepsilon}\}}|f(x,\upsilon_{n})|^{q}dx
+\int_{\{|\upsilon_{n}|>t_{\varepsilon}\}}|f(x,\upsilon_{n})|^{q}dx\\
 &\leq &\displaystyle\omega_{N-1} \max_{B\times [-t_{\varepsilon},t_{\varepsilon}]}|f(x,t)|^{q}+ C\int_{B}e^{\alpha_{0}(\varepsilon+1)|\upsilon_{n}|^{\gamma}}dx.\\
&\leq & C_{\varepsilon}+
 C\displaystyle\int_{B}e^{\alpha_{0}(1+\varepsilon)\|\upsilon_{n}\|^{\gamma}\frac{|\upsilon_{n}|^{\gamma}}{\|\upsilon_{n}\|^{\gamma}}}dx\leq C,
 \end{array}$$
provided $\alpha_{0}(1+\varepsilon)\|\upsilon_{n}\|^{\gamma}\leq p~ \alpha_{N,\beta}$ and
$1<p<\mathcal{U}(u)=(1-\|u\|^{N})^{\frac{-\gamma}{N}}$.\\

\noindent From (\ref{eq:5.17}), there exists $\delta\in (0,\frac{1}{2})$ such that $\rho^{\gamma}=(1-2\delta)\Big(\frac{(\frac{\alpha_{N,\beta}}{\alpha_{0}})^{\frac{N}{\gamma}}}{1-\|u\|^{N}}\Big)^{\frac{1}{(N-1)(1-\beta)}}\cdot$\\
Since $\displaystyle\lim_{n\rightarrow+\infty}\|\upsilon_{n}\|^{\gamma}=\rho^{\gamma}$,
then, for $n$ large enough
 $$\alpha_{0}(1+\varepsilon)\|\upsilon_{n}\|^{\gamma}\leq (1+\varepsilon)(1-\delta)~~\alpha_{N,\beta} \Big(\frac{1}{1-\|u\|^{N}}\Big)^{\frac{\gamma}{N}}.$$
 We choose $\varepsilon>0$ small enough such that $(1+\varepsilon)(1-\delta)< 1$
which implies that
 $$\alpha_{0}(1+\varepsilon)\|\upsilon_{n}\|^{\gamma}<\alpha_{N,\beta} \Big(\frac{1}{1-\|u\|^{N}}\Big)^{\frac{\gamma}{N}} .$$
Hence, the sequence $(f(x,\upsilon_{n}))$ is bounded in $L^{q}$, $q>1$. \\Using the H\"{o}lder inequality, we deduce that
\begin{equation*}\label{eq:5.22}
\begin{array}{rclllll}
\displaystyle \Big|\int_{B}f(x,\upsilon_{n})(\upsilon_{n}-\upsilon)dx\Big| &\leq& \displaystyle\Big(\int_{B}|f(x,\upsilon_{n})|^{q}dx\Big)^{\frac{1}{q}}\Big(\int_{B}|\upsilon_{n}-\upsilon|^{q'}dx\Big)^{\frac{1}{q'}}
\\
 &\leq& \displaystyle C \Big(\int_{B}|\upsilon_{n}-\upsilon|^{q'}dx\Big)^{\frac{1}{q'}}\rightarrow
0~~as~~ n\rightarrow+\infty, \end{array}\end{equation*}
where $~\frac{1}{q}+\frac{1}{q'}=1.$
\\ Since $\langle \mathcal{J_{\lambda}}'(\upsilon_{n}),(\upsilon_{n}-\upsilon)\rangle=o_{n}(1)$,  it follows that
 $$\int_{B}(\omega(x)|\nabla \upsilon_{n}|^{N-2}\nabla \upsilon_{n}.(\nabla \upsilon_{n}-\nabla \upsilon)dx)\rightarrow 0.$$
  On the other side,
  \begin{align*}\int_{B}\omega(x)|\nabla \upsilon_{n}|^{N-2}\nabla \upsilon_{n}.(\nabla \upsilon_{n}-\nabla \upsilon)dx
  =\|\upsilon_{n}\|^{N}-\int_{B}\omega(x)|\nabla \upsilon_{n}|^{N-2}\nabla \upsilon_{n}.\nabla \upsilon dx\cdot\end{align*}
  Passing to the limit in the last equality, we get
  $$\rho^{N}-\|\upsilon\|^{N}=0,$$
therefore $\|\upsilon\|= \rho$. This is in contradiction with (\ref{eq:5.12}). Therefore, $\mathcal{J_{\lambda}}(\upsilon)=c_{\lambda}$.
 By Claim 1, $\mathcal{J_{\lambda}'}(\upsilon)=0$ and by Claim 2, $\upsilon\neq 0$.

\bigskip

\textbf{Proof of Theorem \ref{th1.3}. }From Lemma \ref{lem6} and Lemma \ref{lem10}, we deduce that
$\upsilon$ is a least energy sign-changing solution for problem $(P_{\lambda})$ with exactly tow nodal domains.

\section*{ Declaration of competing interest}
the authors declare that they have no known competing financial interests or personal relationships that could have appeared to influence the work reported in this paper.

\section*{Availability of data}
Data openly available in a public repository that issues data sets with DOIs. We also mention that the documentation to support this study are available from Umm Al-Qura University.

\end{document}